\documentclass[12pt]{amsart}
\usepackage{amssymb}
\usepackage[dvips]{graphics}
\ifx\mathrm\undefined\let\mathrm\rm\fi
\ifx\mathbf\undefined\let\mathbf\bf\fi
\ifx\mathfrak\undefined\let\mathfrak\frak\fi
\ifx\mathcal\undefined\let\mathcal\cal\fi
\ifx\mathbb\undefined\let\mathbb\Bbb\fi
\ifx\emph\undefined\let\emph\it\fi
\font\bb=msbm10 at9.98pt
\begin{document}
\def\semidirect{\hbox{$\;$\bb\char'156$\;$}}
\newcommand{\SL}{\mathrm{SL}}
\newcommand{\ls}{\mathrm{sl}}
\newcommand{\GL}{\mathrm{GL}}
\newcommand{\g}{{{\mathfrak g}\,}}
\newcommand{\bor}{{{\mathfrak b}}}
\newcommand{\n}{{{\mathfrak n}}}
\newcommand{\h}{{{\mathfrak h\,}}}
\newcommand{\Id}{{\operatorname{Id}}}
\newcommand{\Z}{{\mathbb Z}}
\newcommand{\ZZ}{{\mathbb Z_{>0}}}
\newcommand{\N}{{\mathbb N}}
\newcommand{\R}{{\mathbb R}}
\newcommand{\p}{{\mathbb P}}
\newcommand{\C}{{\mathbb C}}
\newcommand{\Q}{{\mathbb Q}}
\newcommand{\D}{\mathcal{D}}
\newcommand{\I}{\mathcal{I}}
\newcommand{\LL}{\mathcal{L}}
\newcommand{\F}{\mathcal{F}}
\newcommand{\W}{\mathcal{W}}
\newcommand{\PP}{\mathcal{P}}
\newcommand{\Sym}{{\rm Sym}}
\newcommand{\Sing}{{\rm Sing}}
\newcommand{\Poly}{{\rm Poly}}
\newcommand{\Span}{{\rm Span}}
\newcommand{\Res}{{\rm Res}}
\newcommand{\1}{{\bf 1}}
\newcommand{\A}{{\bf a}}
\newcommand{\m}{{\bf m}}
\newcommand{\z}{{\bf z}}
\newcommand{\uu}{{\bf u}}
\newcommand{\vv}{{\bf v}}
\newcommand{\w}{{\bf w}}
\newcommand{\kk}{{\bf k}}
\newcommand{\T}{{\bf t}}
\newcommand{\dd}{{\bf d}}
\newcommand{\dontprint}[1]
{\relax}
\newtheorem%
{thm}{Theorem}
\newtheorem%
{prop}
{Proposition}
\newtheorem%
{lemma}
{Lemma}
\newtheorem%
{lemmadef}[thm]{Lemma-Definition}
\newtheorem%
{cor}
{Corollary}
\newtheorem%
{conj}
{Conjecture}
\newenvironment{definition}
{\noindent{\bf Definition\/}:}{\par\medskip}

\title {A theorem of Heine--Stieltjes,\\ the Schubert calculus,\\
and Bethe vectors in the  $\ls_p$ Gaudin model}

\author[{}]
{I.~Scherbak}

\begin{abstract}${}$
Heine and Stieltjes in their studies of linear second-order differential equations with 
polynomial coefficients  having a polynomial solution of a preassigned degree, 
discovered that the roots of such a solution are the coordinates of 
a critical point of a certain remarkable symmetric function,
\cite{He}, \cite{St}. Their result can be reformulated in terms of the 
Schubert calculus as follows:  the critical points label the elements of the intersection 
of certain  Schubert varieties  in the Grassmannian of  two-dimensional subspaces 
of the space of complex polynomials,   \cite{S1}.

In a hundred years after the works of Heine and Stieltjes, it was established 
that the same critical points  determine the Bethe vectors in the $\ls_2$ 
Gaudin model, \cite{G}.
Recently it was proved that the Bethe  vectors of the $\ls_2$ Gaudin model
form a basis  of the subspace of singular vectors of  a given weight in the
tensor product of irreducible $\ls_2$-representations, \cite{SV}.

In the present work we generalize the result of  Heine and Stieltjes to
linear differential equations of order $p>2$. The function, which 
determines elements in the intersection of  corresponding Schubert varieties
 in the Grassmannian of  $p$-dimensional subspaces, turns out to be 
 the very function which appears  in the  $\ls_p$ Gaudin model. 

In the case when the space of states of  the Gaudin model is the tensor product 
of symmetric powers of the standard $\ls_p$-representation, 
we prove that the Bethe  vectors form a basis of the subspace 
of singular vectors of  a given weight. 
\end{abstract}

\maketitle
\pagestyle{myheadings}
\markboth{I.~Scherbak}
{A theorem of Heine--Stieltjes, the Schubert calculus, 
and the $\ls_p$ Gaudin model }
\medskip

\section{Introduction}\label{s1}
The theory of $p$-th order Fuchsian differential equations with polynomial solutions
is closely related to the Schubert calculus in the Grassmannian of $p$-dimensional
subspaces of the vector space of complex polynomials. We describe this relation,
which  is crucial in  the present work,  in Sec.~\ref{s11}.  Then in Sec.~\ref{s12} we 
formulate  a classical result of Heine and Stieltjes in terms of the Schubert calculus, 
 as well as its strengthening  obtained in  \cite{SV}, \cite{S1}. This corresponds to the case $p=2$.
 In Sec.~\ref{s13} we formulate our extension of the Theorem of Heine--Stieltjes to $p>2$. 
 Sec.~\ref{s14} contains an application of  our result to the Gaudin model, 
and  Sec.~\ref{s15} is devoted to the case of special Schubert intersections
in which the result can be strengthened.   

\subsection{Linear differential equations with  polynomial
 solutions and the Schubert calculus.}\label{s11}
See \cite{F} for the Schubert calculus and \cite{R} for the theory of Fuchsian equations.

 If all solutions to a linear differential equation  of order $p$ with meromorphic coefficients
 are polynomials, then  the solution space $V$ is an element 
 of  the Grassmannian  $G_p({\Poly_d})$ of  $p$-dimensional subspaces 
in the vector space ${\Poly}_d$ of  complex polynomials of degree at most $d$, 
where $d$ is the degree of the generic solution. Conversely, any $V\in G_p({\Poly_d})$
defines  a linear differential equation $E_V$ with the solution space $V$. 

For  any $z\in\C$, consider the  flag  $\F_\bullet(z)=
\{\F_0(z) \subset \F_1(z)\subset \dots\subset \F_{d}(z)= {\Poly}_{d}\}$,
where $\F_i(z)$  consists of the polynomials  of the form
$a_i(x-z)^{d-i}+\dots +a_{0}(x-z)^{d}$.
Define $\F_\bullet(\infty)$ by  $\F_i(\infty)={\Poly}_i\,$, $0\leq i\leq d$.
 
 If  integers  $w_1,\dots, w_p$ satisfy  
$ d+1-p\geq w_1\geq \dots w_p\geq 0\,$, we call 
$\w=(w_1,\dots,w_p)$ {\it a Schubert index}.
 For any  $\tilde z\in \C\cup\infty$,   the {\it Schubert cell}
 $\Omega^\circ_{\w}(\tilde z)\subset G_p({\Poly})_d$
 is formed  by  elements  $V\in  G_p({\Poly})_d$ which satisfy
the conditions
 $$
\dim \left(V\cap \F_{d-p+i-w_i}(\tilde z)\right)= i\,,\quad
\dim \left(V\cap \F_{d-p+i-w_i-1}(\tilde z)\right)= i-1\,,
\quad i=1,\dots, p\,.
$$
 The {\it Schubert variety} 
 $\Omega_{\w}(\tilde z)$ is the closure of the Schubert cell,
$$
\Omega_{\w}(\tilde z)=\left\{\, V\in G_p({\Poly}_d)\,
\vert\, \dim \left(V\cap \F_{d-p+i-w_i}(\tilde z)\right)\geq i\,,
\quad i=1,\dots, p\, \right\}\,.
$$
The complex codimension of  $\Omega_{\w}(\tilde z)$ in 
$G_p({\Poly}_{d})$ is
$|\w|=w_1+\dots +w_p\,$.  The homology class $[\Omega_{\w}]$  
of $\Omega_{\w}(\tilde z)$ does not depend on the choice of  flag.
We will denote  $\sigma_{\w}$  the  corresponding  {\it Schubert class},
that is the cohomology class in
 $H^{2|\w|}(G_p({\Poly}_d))$, whose cap product with the
fundamental class of   $G_p({\Poly}_d)$ is  $[\Omega_{\w}]$, see \cite{F}.

For any  $\tilde z\in \C\cup\infty$ and any $V\in G_p({\Poly}_{d})$, 
 we have $V\in \Omega^\circ_{\w(\tilde z;V)}(\tilde z)$ for a certain 
Schubert index $\w(\tilde z;V)$. Equivalently,  the exponents of  the equation $E_V$ 
at $\tilde z$  are   $w_{p+1-i}(\tilde z;V)+i-1$,  $1\leq i\leq p\,$, see Sec.~\ref{s23}.  

{\it The Wronskian}  of  $V\in G_p({\Poly_d})$ is defined as a monic polynomial 
$W_V(x)$ which is proportional to the Wronskian of   some (and hence, any) 
basis of $V$. See \cite{EGa}, \cite{KhSo} for relations of  the Wronskian to the Schubert
calculus. Clearly $W_V(x)$ is also the Wronskian of the differential  equation $E_V$.

\medskip
If  $W_V(z_0)\neq 0$, then  $\w(z_0;V)=(0,\dots,0)$.
Equivalently, any singular point of the  
equation $E_V$ is a root of the Wronskian $W_V(x)$.
Moreover, for any singular point $\tilde z$
the codimension of   $\Omega_{\w(\tilde z;V)}(\tilde z)$ 
is exactly the multiplicity of  $\tilde z$  as a root of the Wronskian.
The sum of  the codimensions of  the Schubert varieties corresponding
to all singular points (including $\infty$)  equals the dimension of  
$G_p({\Poly}_{d})$  (this is in fact the Fuchs theorem on exponents 
at singular points, see  \cite{R}). One can see that the intersection 
of these  varieties coincides with the intersection of the corresponding 
Schubert cells.

\subsection {Theorem of Heine--Stieltjes in terms of the 
Schubert calculus.}\label{s12}
In  \cite{He}, pp.~472--479, \cite{St}, 
Heine and Stieltjes studied the following problem. 

\medskip\noindent
 {\it PROBLEM\ \ Let $A(x)$ and $B(x)$ be given complex polynomials
of degrees $n$ and $n-1$ respectively, $n\geq 2$. To determine a polynomial 
$C(x)$ of degree less than or equal to $n-2$ such that the equation
$$A(x)u''(x)+B(x)u'(x)+C(x)u(x)=0$$
has a solution which is a complex polynomial of a preassigned
degree $k$}.

\medskip
They assumed that  $A(x)$ has no multiple roots,
 $A(x)=\prod_{j=1}^{n}(x-z_j)$. Then  $B(x)$ can be given by 
${B(x)}/{A(x)}=-\sum_{j=1}^{n} {m_j}/{(x-z_j)}$, where $m_j$ 
and $z_j$ are suitable complex numbers. In the present work
we are interested in the result of Heine and Stieltjes under the following
additional assumption,

\medskip\centerline
{\it all  numbers  $m_1,\dots, m_n$ are positive integers.}

\medskip\noindent
In this case the theory of Fuchsian equations (\cite{R}, Ch.~6)  asserts  that if  the equation  
has a polynomial solution $u(x)$ with no multiple roots, then {\it all}  solutions to 
this equation are polynomials and $u(z_j)\neq 0$, $1\leq j\leq n$.

 In order to formulate the result of   Heine and Stieltjes in terms of the Schubert calculus, 
we call  $V\in G_2({\Poly}_{d})$ {\it a nondegenerate two-plane} if  its 
polynomials of the smallest degree do not have multiple roots.  
Write $m=(m_1,\dots,m_n)$, $M=m_1+\dots +m_n$,  $z=(z_1,\dots,z_n)$ and 
\begin{equation}\label{W}
W_{m,z}(x)=(x-z_1)^{m_1}\dots (x-z_n)^{m_n},\quad  \deg W_{m,z}=M\,.
\end{equation}
If $W_V(x)=W_{m,z}(x)$ and if $V$  contains a polynomial of  degree $k$, 
then $V$ contains polynomials  of degree  $M +1-k$ as well.
Denote  $\Delta(u)$  the discriminant of the polynomial $u(x)$ and   
${\rm Res}(u,W_{m,z})$  the resultant of $u(x)$ and  $W_{m,z}(x)$. 
If $u(x)$ is a polynomial of the smallest degree in a nondegenerate two-plane,
then $\Delta(u)\neq 0.$

\medskip
{\bf  Theorem ~A} (Heine--Stieltjes)\  \ {\it Let $u(x)$ be an unknown polynomial 
of positive degree $k\leq M/2$. 
Then the critical points with non-zero critical values of the function
$$\Phi_{m,z}(u)=\frac{\Delta (u)}{\Res(u,W_{m,z})}$$ 
are in a one-to-one correspondence
with the nondegenerate two-planes in the  
intersection of Schubert varieties
$$
\I_m(z)=
\Omega_{(m_1,0)}(z_1)\cap\dots\cap 
\Omega_{(m_n,0)}(z_n)\cap
\Omega_{(M-2k;0)}(\infty)
\subset G_2({\Poly}_{M+1-k})\,,
$$
and  determine solutions to the Problem.}

\medskip
Intersections of  Schubert varieties corresponding to  flags $\F_\bullet(\tilde z)$
were studied by D.~Eisenbud and J.~Harris. 
The results of  \cite{EH} say that   $\I_m(z)$
is zero-dimensional and  consists of  
$\sigma_{(m_1,0)}\cdot \, ...\, \cdot\sigma_{(m_1,0)}\cdot \sigma_{(M-2k,0)}$
elements counted with multiplicities.  

\medskip
If $z$ is generic, then one can say more.
Here and  everywhere in the text, the words ``generic $z$''  mean   
that $z$ does not belong  to a suitable proper algebraic  surface in $\C^n$.

\medskip
{\bf  Theorem ~B} (\cite{S1}, cf. Theorem 12 of  \cite{SV})\ \ 
{\it For generic $z$, the intersection $\I_m(z)$  is transversal and
consists of   nondegenerate two-planes only.}

\medskip
Thus the  critical points of  the function $\Phi_{m,z}$ 
 determine the  intersection  $\I_m(z)$ and  {\it all} solutions to the Problem,
 for generic $z$.
The number of  these critical points was calculated explicitly in  \cite{SV}, 
Theorems~1 and~5.

\subsection {Nondegenerate $p$-planes}\label{s13}
In the present work we  extend  Theorems A and B to $p>2$.
Fix    $ \{\w\}=\{\w(1), \dots ,\w(n),  \w(n+1)\}$, the set of $n+1$ Schubert indices  
$\w(j)=(w_1(j),\dots w_{p-1}(j), 0)$, such that
$$|\w(j)| = m_j\,, 1\leq j\leq n,\quad |\w(n+1)|=\dim  G_p({\Poly}_{d})-\sum_{j=1}^nm_j.$$
Consider  the intersection of  Schubert varieties
\begin{equation}\label{I}
\I_{\{\w\}}(z)=\Omega_{\w(1)}(z_1)\cap \dots
\cap \Omega_{\w(n)}(z_n)\cap
\Omega_{\w(n+1)}(\infty)\subset G_p({\Poly}_d)\,.
\end{equation}
This intersection is zero-dimensional, \cite{EH}. Moreover,
if $V\in \I_{\{\w\}}(z)$, then $V$ has no base point and 
its Wronskian is $W_{m,z}(x)$ given by  (\ref{W}).

For any  $V\in \I_{\{\w\}}(z)$  denote   
 $V_\bullet=\{ V_1\subset V_2 \subset \dots \subset V_p=V\}$
the  flag  obtained by the  intersection of  $V$
and  $\F_\bullet(\infty)$. Denote $W_i(x)$
the Wronskian of $V_i$, $1\leq i\leq p$.
In particular, $W_p(x)=W_{m,z}(x)$. We say that $V$ is {\it a nondegenerate $p$-plane}
in $\I_{\{\w\}}(z)$ if 
\begin{itemize}
\item 
for any $\tilde z\in\{z_1,\dots, z_n, \infty\}$, any $1\leq i\leq p$, 
and any $1\leq l\leq i$, there exists a polynomial in $V_i$  
which has a root of  order {\it exactly} $w_{p+1-l}(\tilde z;V)+l-1$  at $\tilde z$;
\item
for any  complex number $t\notin \{z_1,\dots, z_n\}$ and any  
$1\leq i\leq p-2$ such that
$W_i(t)=0$, we have $W_{i+1}(t)\neq 0$.
\end{itemize}

In order to write down the  function whose critical points determine 
the nondegenerate  $p$-planes  in $\I_{\{\w\}}(z)$,  
introduce  notions of   {\it  relative discriminant} and {\it relative resultant}.
For fixed $z=(z_1,\dots,z_n)$,  any monic polynomial 
$P(x)$ can be presented 
in a unique way as the product of two monic polynomials $T(x)$ and $Z(x)$ which
satisfy 
$$
P(x)=T(x)Z(x), \ \   T(z_j)\neq 0,\ \  Z(x)\neq 0  
{\rm\ \ for\ any\ } x\neq z_j, \ \ 1\leq j\leq n.
$$
Define {\it  the relative discriminant of $P(x)$  with respect to $z$} as
$$\Delta_{z}(P)=\frac{\Delta(P)}{\Delta(Z)}=\Delta(T){\rm Res}^2(Z,T),$$
and {\it the relative resultant   of  $P_i(x)=T_i(x)Z_i(x)$, $i=1,2$,
with respect to $z$} as
$${\rm Res}_{z}(P_1, P_2)=\frac{{\rm Res}(P_1,P_2)}{{\rm Res}(Z_1,Z_2)}=
{\rm Res}(T_1, T_2){\rm Res}(T_1, Z_2){\rm Res}(T_2, Z_1).$$
For  $V\in \I_{\{\w\}}(z)$,
take the flag $V_\bullet$ and present the Wronskian $W_i(x)$ of $V_i$ in such a form, 
$$
W_i(x)=Z_i(x)T_{p-i}(x).
$$
If  $V$ is a  nondegenerate  $p$-plane, then the  Schubert indices $\{\w\}$ define 
polynomials $Z_1,\dots Z_p$  and the degrees $k_0,\dots, k_{p-1}$
of  polynomials $T_0,\dots , T_{p-1}$, see  Sec.~\ref{s31}.
In particular,  $Z_p(x)=W_{m,z}(x)$,  $Z_1(x)=1$, $T_0(x)=1$, i.e. $k_0=0$.
Moreover, polynomials $T_1, \dots, T_{p-1}$ do not have multiple roots and
any two neighboring of them do not have common roots.

\medskip {\bf Definition}\ \ The function
\begin{equation}\label{Phi-w}
\Phi_{ \{\w\},z} (T_1,\dots, T_{p-1})=
\frac {\Delta_z(W_1)\cdots \Delta_z(W_{p-1})}
{{\rm Res}_z(W_1,W_2)\cdots
{\rm Res}_z(W_{p-1},W_{p})}
\end{equation}
is  called {\it the generating function of the Schubert intersection} 
 $\I_{\{\w\}}(z)$.

\medskip
{\bf Theorem 1}\ \ {\it There is a one-to-one correspondence between
the critical points with non-zero critical values of the function 
$\Phi_{ \{\w\},z} (T_1,\dots, T_{p-1})$
and the nondegenerate $p$-planes in}  $\I_{\{\w\}}(z)$. 

\medskip
The proof of  Theorem~1 is given in Sec.~\ref{s33}.

\medskip
{\bf  Remark}\ \ For {\it any}  $V\in \I_{\{\w\}}(z)$,   any
$1\leq i\leq p-1$,  and any $\tilde z\in\{z_1,\dots, z_n, \infty\}$, 
there exists a basis $P_(x),\dots,P_i(x)$ in $V_i$ such that $P_l(x)$  
has a root of  order {\it at least} $w_{p+1-l}(\tilde z;V)+l-1$  at $\tilde z$,
$1\leq l\leq i$. This means that the Wronskian 
$W_i(x)=W_{V_i}(x)$ is of the  form $W_i(x)=Z_i(x)F_{p-i}(x)$, where $Z_i(x)$ 
is defined by the Schubert indices as above and $F_{p-i}(x)$ is a certain polynomial.  
The conditions for $V$ to be a nondegenerate $p$-plane mean that 
{\it the polynomials  $F_i(x)$, $1\leq i\leq p-1$, are  as generic as possible}.

\medskip
{\bf Conjecture}\ \ {\it For generic  $z$, the intersection   $\I_{\{\w\}}(z)$
given by {\rm (\ref{I})} is transversal and consists of  nondegenerate $p$-planes only; 
in particular, the critical points  of the  function 
$\Phi_{ \{\w\},z}(T_1,\dots, T_{p-1})$ 
determine   all elements of  \  $\I_{\{\w\}}(z)$, for generic $z$}. 

\medskip
Plausibly,   the statements  
{\it ``$V\in \I_{\{\w\}}(z)$ is a nondegenerate  $p$-plane''}
and  {\it ``$V$ is a simple point of  the Scubert intersection $\I_{\{\w\}}(z)$''}  
are  equivalent; in other words, the conditions of nondegenericity provide
intrinsic characterization for  the points of transversal intersection.

\medskip
We succeed to prove the  Conjecture only in the case  of {\it special}  Schubert 
intersections  using the relation to the Gaudin model; 
these are Theorem~2 and its Corollary formulated in Sec.~\ref{s15}.

\subsection{ Bethe vectors in the $\ls_p$ Gaudin model}\label{s14}
On Bethe vectors in the Gaudin model  see \cite{G}, \cite{FeFrRe},  
\cite{Fr}, \cite{ReV}.

Consider Lie algebra $\ls_p=\ls_p(\C)$.
Fix  the tensor product $L=L_1\otimes\dots\otimes L_n$  of irreducible
finite-dimensional $\ls_p$-representations and $z=(z_1,\dots, z_n)$ with
pairwise distinct complex coordinates. {\it The Gaudin hamiltonians}  
associated with $z$ and  $L$,  are defined as follows,
\begin{equation}
\label{H}
H_i(z)=\sum_{j\neq i}\frac{C_{ij}}{z_i-z_j}\,,
\quad  i=1,\dots, n,
\end{equation} 
where   $C_{ij}\,:\, L\to L$ is  a linear
operator acting  as the Casimir element on the factors
$L_i$ and $L_j$ and as the identity on the other factors of $L$.
The  Gaudin hamiltonians commute, and one of the main  problems
in the Gaudin model is to find their common eigenvectors.
This problem is solved by means of
{\it the Bethe Ansatz} which  is a widely used method for diagonalizing
of  commuting  hamiltonians in integrable models of statistical mechanics.
The idea of the Bethe Ansatz is to construct a certain   function  
$v\,:\,\C^K\,\to\, L$
in such a way that for some special values  of  the  argument
$\T=(t_1,\dots,t_K)$, the values of the function are eigenvectors. 
The  equations which determine
these special values of the argument are called {\it the Bethe equations}
and the corresponding vectors  $v(\T)$ are called {\it the Bethe vectors}. 
The subspace of singular vectors of a given  weight in $L$  is an invariant 
subspace of the Gaudin hamiltonians. The  Bethe equations in the Gaudin
model associates with $z$ and with this subspace  provide the critical point system 
of  a certain function called {\it the master function} of the model, \cite{G},
\cite{FeFrRe}, \cite{ReV}.

\medskip
 It turns out  that the function  (\ref{Phi-w})  written in terms of unknown roots of 
 polynomials $T_1,\dots, T_{p-1}$ is the master function  of the Gaudin model
associated with $z$ and the subspace  $\Sing_{\kk}\Gamma_{\{\A\}}$ 
of singular vectors of the weight 
$$
\Lambda(\kk)=\Lambda_{\A(1)}+\dots + \Lambda_{\A(n)}-
k_1\alpha_1-\dots-k_{p-1}\alpha_{p-1}\,, \ \kk=(k_1,\dots,k_{p-1}), 
\ \ k_i=\deg T_i\,, 
$$                          
in the tensor product 
$\Gamma_{\{\A\}}=\Gamma_{\A(1)}\otimes\dots\otimes \Gamma_{\A(n)}$
of  irreducible $\ls_p$-representations $\Gamma_{\A(j)}$ with integral
dominant highest weights  $\Lambda_{\A(j)}$,  
$$
(\Lambda_\A(j), \alpha_i)=a_i(j),\ \  a_i(j)=w_i(j)-w_{i+1}(j),\ \  
1\leq i\leq p-1,\ \  1\leq j\leq n,
$$
here  $\alpha_1,\dots, \alpha_{p-1}$  are simple roots of the Lie algebra
$\ls_p=\ls_p(\C)$ and $(\cdot,\cdot)$ is the Killing form on the dual to the
Cartan subalgebra.

\medskip
For any $1\leq j\leq n$, 
the Schubert class $\sigma_{\w(j)}$ and the $\ls_p$-module 
$\Gamma_{\A(j)}$ have the same  Young diagram,  \cite{F}.
The condition that $\Lambda(\kk)$ is an integral 
dominant weight is necessary for  $\I_{\{\w\}}(z)$ might be non-empty. 
Under this condition, 
the Schubert class dual to $\sigma_{\w(n+1)}$ and   the $\ls_p$-module 
$\Gamma_{\A(n+1)}$ have the same  Young diagram.
The relation of the Schubert calculus to 
the representation theory of the Lie algebra $\ls_p(\C)$
implies that the dimension of  $\Sing_{\kk}\Gamma_{\{\A\}}$ 
is the intersection number of the Schubert classes
$\sigma_{\w(1)}\cdot {\rm \ ... \ }\sigma_{\w(n)}
\cdot \sigma_{\w(n+1)}$, \cite{F}.  Theorem~1 leads to the following
claim.

\medskip
{\bf  Corollary from Theorem~1}\ \
{\it For the Gaudin model,  which
associates with $z$ and $\Sing_{\kk}\Gamma_{\{\A\}}$, 
the number of Bethe vectors is at most  
$\dim \Sing_{\kk}\Gamma_{\{\A\}}$.}

\medskip
{\bf Remark}\ \  
While the first version of  the text was in preparation,    
E.~Mukhin and A.~Varchenko  have  posted on arXiv  a preprint
containing this result  in a different formulation, see  Sec.~5.6 of \cite{MV}.

\subsection{Special Schubert intersections and the Gaudin model
associated to the tensor product of symmetric powers of the standard
$\ls_p$-representation}\label{s15}
Let the Schubert indices $\w(1), \dots , \w(n)$  be {\it special}, 
$$
\w(j)=(m_j,0,\dots,0)\,,\ 1\leq m_j\leq d+1-p\,,\ 1\leq j\leq n\,.
$$
Denote  $\Omega_{(m_j)}=\Omega_{(m_j,0,\dots, 0)}$  the
 corresponding {\it special} Schubert varieties.  
The intersection given by (\ref{I}) becomes 
\begin{eqnarray}\label{Im}
&&\I_{\{m,\w(n+1)\}}(z)=\Omega_{(m_1)}(z_1)\cap \dots
\cap \Omega_{(m_n)}(z_n)\cap
\Omega_{\w(n+1)}(\infty)\subset G_p({\Poly}_d)\,,\\
&& m=(m_1,\dots, m_n),\   \  |\w(n+1)|=\dim  G_p({\Poly}_{d})-\sum_{j=1}^nm_j\,,\nonumber
\end{eqnarray}
which is a {\it  special} Schubert intersection. 
 
In this case  $Z_i(x)=1$,  and hence $T_{p-i}(x)=W_i(x)$ for all $1\leq i\leq p-1$.
Therefore the relative discriminants  and resultants turn into the usual ones and
the  generating function takes the  form
\begin{equation}\label{G}
\Phi_{\{m, \w(n+1)\}, z}(W_1,\dots,W_{p-1})\,=\,
\frac {\Delta(W_1)\cdots \Delta(W_{p-1})}
{{\rm Res}(W_1,W_2)\cdots
{\rm Res}(W_{p-1},W_{p})}\,.
\end{equation}
This is the master function of the Gaudin model associated to 
the tensor product of the $m_j$-th symmetric powers of the standard 
$\ls_p$-representation,  and the following result holds.

\medskip
{\bf  Theorem~2} \ \ {\it In the Gaudin model associated with 
 the tensor product of symmetric powers of the standard $\ls_p$-representation,
 the Bethe vectors form a basis in the subspace of singular vectors of a given
 weight, for generic $z=(z_1,\dots,z_n)$.}

\medskip
The proof is given in Sec.~\ref{s53}. This theorem immediately implies
the statement of the Conjecture for the special  Schubert intersections.

\medskip
{\bf   Corollary from  Theorem~2}\ \ {\it For generic $z$, the  special Schubert
intersection {\rm (\ref{Im})} is transversal and consists of nondegenerate 
$p$-planes only.}

\medskip 
Since the   intersection $\I_m(z)$ of Theorem~B is special,  
Theorem~2 and its Corollary can be considered as an extension of  Theorem~B.

\bigskip\noindent
{\bf  Plan of the paper}\ \

Sec.~\ref{s2} is an exposition of the theory of linear differential equations with
only polynomial solutions in terms of the Schubert calculus and Wronskians.
Results of this section are Proposition~\ref{PC1}  of Sec.~\ref{s23} containing
an upper bound for the  number of the equations with given singular points
and exponents, and Proposition~\ref{P3} of Sec.~\ref{s24} which is an essential
ingredient in the proof of Theorem~1 given in Sec.~\ref{s3}. Corollary~\ref{C11}
of Sec.~\ref{s34} is a reformulation of Theorem~1 in terms of rational curves
in $\C\p^{p-1}$.  In Sec.~\ref{s4} we collect necessary  facts on the 
$\ls_p$-representations and on the $\ls_p$ Gaudin model and deduce the 
Corollary from Theorem~1. In Sec.~\ref{s43} we discuss the problem on 
the simplicity of the common spectrum of the Gaudin hamiltonians. 
Sec.~\ref{s5} is devoted to the case of special Schubert intersections. 
Propositions~\ref{P7} and~\ref{PC5} of Sec.~\ref{s51} give an explicit form of
the corresponding differential equations. Using these propositions and relations
to the Gaudin model, we prove Theorem~2 and deduce Corollary~\ref{C9}  on
the number of the corresponding differential equations.

\bigskip\noindent
{\bf Acknowledgments}\ \ 

I am  grateful to A.~Gabrielov for sharing with me and allowing to use
his unpublished note on  the Wronski map contained  formula ~(\ref{G}) and 
an  algebraic proof of the equality~(\ref{Ga})  of Proposition \ref{P3}.
 
It is my pleasure to  thank E.~Frenkel, D.~Kazhdan and  F.~Sottile for 
interesting and fruitful discussions.

Finally, I  would like to take the opportunity to thank A.~Varchenko.  
The present work can be considered as a continuation of our joint work~\cite{SV} 
where the case $p=2$, with no relation to the Schubert calculus, was studied; 
another independent continuation of~\cite{SV}, dealing with Kac--Moody algebras 
and flag varieties, is contained in~\cite{MV}.
 
\section{Fuchsian equations with only polynomial solutions\\
and the Schubert calculus}\label{s2}
\subsection{Linear differential equations
with only polynomial solutions}\label{s21}

Denote ${\Poly}$ the vector space of  complex polynomials in one
variable and $G_p({\Poly})$  the Grassmannian of  $p$-dimensional
subspaces of ${\Poly}$. Fix $V\in G_p({\Poly})$ and consider
the linear differential equation  $E_V$ with meromorphic coefficients
and the solution space $V$.
This equation  is unique, up to multiplication by
a meromorphic function. Let  $u=u(x)$ be an unknown
function.
If  polynomials  $P_1(x),\dots, P_p(x)$ form a basis of $V$, then
$E_V$ can be written in the form
\begin{equation}\label{E1}
\det \left(\begin{array}{cccc} u(x)\ \ &P_1(x)\   \ &\dots\   \ &P_p(x)\ \ \\
                          u'(x)\ \ &P_1'(x)\ \  &\dots\ \ &P_p'(x)\    \\
                           \dots\ \  &\dots\ \   &\dots\ \  &\dots\ \\
                          u^{(p)}(x)\   &P_1^{(p)}(x)\ \ &\dots\ \ &P_p^{(p)}(x)\\
                           \end{array}\right) = 0\,.
\end{equation}
On the left hand side  the Wronski determinant of  functions  
$u(x)$,  $P_1(x),\dots,P_p(x)$ stands.

\medskip
The equation $E_V$ is {\it Fuchsian}.  On Fuchsian equations
see, for example, Ch.~6 of~\cite{R}. Any  solution of  any  Fuchsian
equation at any point $z\in\C$ has the form
$$
u(x)=(x-z)^{\rho}\sum_{l=0}^{\infty}c_l(x-z_0)^l,\quad c_0\neq 0,
$$
where $\rho$ is a suitable complex number called {\it an exponent at} $z$.
A point where all coefficients of the equation are holomorphic
is {\it an ordinary point}. At any ordinary point the exponents
are $0,1,\dots, p-1$. A non-ordinary point  is called {\it singular}.

A Fuchsian equation has {\it an ordinary (resp., singular)
point at infinity} if after the change
$x=1/\xi$ the point $\xi=0$ becomes an ordinary (resp., singular)
point of the transformed equation.

\medskip
For  any $z\in\C$,   an exponent of the equation $E_V$
at $z$  is a non-negative integer $\rho=\rho(z)$ such that the
equation has a solution of the form $(x-z)^{\rho}f(x)$, where
$f(x)$ is a polynomial which does not vanish at $z$.
An  exponent at infinity is a non-positive integer $-k$
such that  the equation has a polynomial solution of degree $k$.
For  any $\tilde z\in\C\cup\infty$,
the set of exponents at $\tilde z$ consists of $p$ pairwise distinct integers.

\medskip
Any finite singular point of the equation  $E_V$
is a root of the Wronskian $W_V(x)$.
If $z_0\in\C$ is a root of $W_V(x)$ of multiplicity $m_0$,
then the sum of all exponents at $z_0$ is $m_0+p(p-1)/2$.

\medskip
Assume that the generic solution to the  equation $E_V$
has degree $d$. This means that $V$ belongs to
$G_p({\Poly}_d)$, where ${\Poly}_d\subset{\Poly}$
is the subspace of polynomials of degree at most $d$.
In this case the degree of  $W_V(x)$ is at most  $p(d+1-p)$. We say
that $p(d+1-p)-\deg W_V$ is {\it the multiplicity  of the root
at $\infty$}.  The sum of all exponents at infinity  is $-\deg W_V
-p(p-1)/2$.
If the multiplicity at $\infty$ is $0$, then the equation
does not have singularity at $\infty$ and the exponents at $\infty$
are $-d+p-1,-d+p-2,\dots,-d$.

\medskip
Any Fuchsian differential equation of order $p$ with only polynomial
solutions is defined by its solution space.  For $V\in G_p({\Poly})$,
call the singular points of $E_V$ {\it the singular points of $V$} and
the exponents of  $E_V$ at
$\tilde z\in\C\cup\infty$ {\it the exponents of $V$ at $\tilde z$}.

\medskip\noindent
{\bf Remark}\ \
All results related to  Fuchsian differential equations with only
polynomial solutions can be reformulated for  Fuchsian  equations
with only univalued solutions,  considered up to the rational changes.
Indeed, consider  a Fuchsian equation of order $p$
with respect to unknown function $u$ with only univalued, that is  rational,
solutions.  The solution space can be spanned by $p$ functions of the form
$P_i/Q$, where $P_1,\dots, P_p$ and $Q$ are polynomials.
The change of dependent variable $u\mapsto u/Q$ gives a new  Fuchsian
equation,  and the solution space of the new equation is spanned by polynomials 
$P_1,\dots, P_p$. Cp.  Sec.~1.3 of  \cite{SV}.

\subsection{Intersections of Schubert classes}\label{s22} (See, for example, \cite{F}.)
We use notation of Sec.~\ref{s11}. 
The Grassmannian  $G_p({\Poly}_d)$ is a smooth complex algebraic variety
of dimension $p(d+1-p)$, and 
the Schubert classes  $\{\sigma_{\w}\}$ for all possible Schubert indices
give a basis over $\Z$ for the cohomology ring $H^*(G_p({\Poly}_d))$.
The product or {\it intersection} of   two Schubert classes
$\sigma_{\w}$ and  $\sigma_{\vv}$
has the form
$$
\sigma_{\w}\,\cdot\,\sigma_{\vv}\,=\,
\sum_{\uu:\, |\uu|=|\w|+|\vv|} C(\w;\vv;\uu)
\sigma_{\uu}\,,
$$
where coefficients $C(\w;\vv;\uu)$ are nonnegative integers
known as {\it the Littlewood--Richardson coefficients}.

\medskip
If the sum of  codimensions of Schubert classes equals the dimension
of $G_p({\Poly}_d)$, then their intersection is an integer,
identifying the generator of the top cohomology group
$$\sigma_{(d+1-p,\dots, d+1-p)}\in H^{2p(d+1-p)}(G_p({\Poly}_d))$$
with $1\in\Z$. This integer is  called {\it the intersection number}.

\subsection{Fuchsian equations with only polynomial solutions
and Schubert intersections  in $G_p({\Poly}_d)$}\label{s23}
Fix $V\in G_p({\Poly}_d)$ and consider the equation $E_V$.
The degree of  the generic solution to this equation is less than or equal
to $d$.
For any  $z\in\C$, write the exponents in the decreasing order,
$d\geq \rho_1(z)>\dots >\rho_p(z)\geq 0$.
Write the  exponents  at $\infty$ in the decreasing order,
 $0\geq -d_1> \dots >-d_p\geq -d$.  Define {\it the Schubert index
 $\w(\tilde z;V)=\left(w_1(\tilde z),\dots,w_p(\tilde z)\right)$
 of $V$ at $\tilde z\in \C\cup\infty$}  as follows,
\begin{equation}\label{E3}
 w_i(z)=\rho_i(z)+i-p\,\quad  {\rm   if\ } z\in\C\,;
\quad   w_i(\infty)=d-d_i+i-p\,, \quad  i=1,\dots,p.
\end{equation}
The theory of Fuchsian equations with only polynomial
solutions of Sec.~\ref{s21}  can be reformulated in terms of the Schubert calculus.

\begin{prop}\label{P1} Let  $V\in G_p({\Poly}_d)$. Then
\begin{itemize}
\item
$V\in\cap_{\tilde z\in\C\cup\infty} 
\Omega^\circ_{\w(\tilde z;V)}(\tilde z) =
\cap_{\tilde z\in\C\cup\infty}\Omega_{\w(\tilde z;V)}(\tilde z)\,;$\\
\
\item
$\tilde z\in\C\cup\infty$ is an ordinary point  of the equation
$E_V$ if and only if  $\w(\tilde z;V)=(0,\dots,0)$;\\
\
\item
for any   singular point $\tilde z\in\C\cup\infty$ of the equation  $E_V$, 
the codimension $|\w(\tilde z;V)|$ is the multiplicity  of $\tilde z$ as 
a root of the Wronskian $W_V(x)$; \\
\
\item
$\sum_{\tilde z\in\C\cup\infty}|\w(\tilde z;V)|\,=\,\dim G_p({\Poly}_d)\,.$
\hfill   $\triangleleft$
\end{itemize}
\end{prop} 

Thus the intersection of all Schubert  varieties $\Omega_{\w(\tilde z;V)}(\tilde z)$ 
coincide with  the intersection  
\begin{equation}\label{IV}
\I_V=\Omega_{\w(z_1;V)}(z_1)\cap\dots\cap
\Omega_{\w(z_n;V)}(z_n)\cap
\Omega_{\w(\infty;V)}(\infty)\,,
\end{equation}
where $z_1,\dots, z_n$ are all finite singular points of $E_V$.
As it is known   (\cite{EH}),  the  intersection  (\ref{IV})  is zero-dimensional
and  its cardinality is bounded from above by the intersection number 
of the corresponding Schubert classes, $\sigma_{\w(z_1;V)}\cdot\, ... \, \cdot
\sigma_{\w(z_n;V)}\cdot\sigma_{\w(\infty;V)}\,$.

If  $V'\in \I_{V}$, then  the equations 
$E_V$ and $E_{V'}$ have the same  singular points and the same exponents 
at these points. We arrive at the following statement on Fuchsian equations.

\begin{prop}\label{PC1}
For any  $V \in G_p({\Poly}_d)$, the  Fuchsian differential equations of order 
$p$ having  the same singular points and the same  exponents at these points as 
$E_V$ are in a one-to-one correspondence with the elements of  the Schubert
intersection $\I_V$ given by {\rm (\ref{IV})}. 
The number of such equations is finite and  bounded from
above by the intersection number of the corresponding Schubert classes.
\hfill   $\triangleleft$
\end{prop}

\subsection{Wronskians}\label{s24}

Fix $V\in G_p({\Poly}_d)$  and a basis  $P_1(x),\dots, P_p(x)$ of $V$.
For $1\leq l\leq p$, denote  $V(l)=\Span\{P_1,\dots,P_l\}$ 
and  $W_l(x)$ the Wronskian of  $V(l)$, $1\leq l\leq p$. In particular,
we have $W_1(x)=P_1(x)$, $V(p)=V$ and $W_p(x)=W_V(x)$.
For convenience, set $W_0(x)=1$.

\begin{prop}\label{P2}{\rm (\cite{PSz}, Part VII, sec.~5, Problem 62)}
The equation $E_V$ can be written in the form
\begin{equation}\label{E4}
\frac{d}{dx}\,\frac{W^2_{p-1}(x)}{W_{p-2}(x)W_p(x)}\, \dots \,
\frac{d}{dx}\,\frac{W^2_2(x)}{W_3(x)W_1(x)}\,\frac{d}{dx}\,
\frac{W^2_1(x)}{W_2(x)W_0(x)}\,
\frac{d}{dx}\,\frac{u(x)}{W_1(x)}=0. \quad    \triangleleft
\end{equation}
\end{prop}

\begin{cor}\label{C3}
For the equation {\rm (\ref{E4})}, the functions
\begin{eqnarray*}
 u_1(x) & = & W_1(x), \\
 u_2(x)&=&W_1(x)\int^x\frac{W_2W_0}{W_1^2}\,,\\
 u_3(x)&=&W_1(x)\int^x\left(\frac{W_2(\xi)W_0(\xi)}{W_1^2(\xi)}
 \int^{\xi}\frac{W_1W_3}{W_2^2}\right)\,,\\
 \dots & \ & \dots \quad \dots\quad \dots  \\
u_p(x) & = & W_1(x)\int^x\left(\frac{W_2(\xi)W_0(\xi)}{W_1^2(\xi)}\int^{\xi}\left(
\frac{W_1(\tau)W_3(\tau)}{W_2^2(\tau)}\int^{\tau}\dots\int^{\eta}\left.
\frac{W_{p-2}W_p}{W^2_{p-1}}\right.\right)\dots\right)
\end{eqnarray*}
form a basis in the solution space, and the polynomial $W_i$ is the Wronskian of
$u_1,\dots, u_i\,,\ \ 1\leq i\leq p$. \hfill   $\triangleleft$
\end{cor}

\begin{prop}\label{P3}
If  for some $t\in\C$ and for some  $1\leq l\leq p-1$,
we have   $W_{l}(t)=0$ and $W_{l\pm 1}(t)\neq 0$, then $W'_l(t)\neq 0$ and
\begin{equation}\label{Ga}
\frac{W''_l(t)}{W'_l(t)}-\frac{W'_{l+1}(t)}{W_{l+1}(t)}-
\frac{W'_{l-1}(t)}{W_{l-1}(t)}=0\,.
\end{equation}
\end{prop}

\noindent{\bf Proof:}\ \  The equation $E_{V(l)}$  can be written as follows,

$$
\frac{d}{dx}\,\frac{W^2_{l-1}(x)}{W_{l-2}(x)W_l(x)}\, \dots \,
\frac{d}{dx}\,\frac{W^2_1(x)}{W_0(x)W_2(x)}\,
\frac{d}{dx}\,\frac{u(x)}{W_1(x)}=0,\ \ l=2,\dots, p.
$$
Any solution to the equation $E_{V(l)}$ is clearly a solution to the
equation $E_{V(l+1)}$ as well.
Proposition~\ref{P1} asserts that
$t$ is an ordinary point of  the equations  $E_{V(l\pm 1)}$  and
a singular point of the equation  $E_{V(l)}$.
Thus the set of exponents at $t$ of the equation
$E_{V(l)}$ contains $\{0,1,\dots,l-2\}$,
is contained in the set $\{0,1,\dots,l\}$,
and differs from the set $\{0,1,\dots,l-1\}$. Hence this set
is $\{0,1,\dots,l-2, l\}$, the Schubert index of $V(l)$ at $t$ is
$\w(t;V(l))=(1,0,\dots,0)$, $\vert\w(t;V(l))\vert=1$, and 
$t$ is a simple root of the Wronskian $W_l(x)$.

In order to prove the equality (\ref{Ga}), note that the solution $u_{l+1}(x)$ 
of the equation $E_{V(l+1)}$ given in Corollary~\ref{C3} is a polynomial. 
Therefore  the most interior integral should be a rational function.
In particular, the residue of the integrand,
$W_{l+1}(x)W_{l-1}(x)/W_l^2(x)$, at $t$ should be zero. We have
\begin{eqnarray*}
&\ & \frac {W_{l+1}(x)W_{l-1}(x)}{W_l^2(x)}= \\
& = &\frac{\left(W_{l+1}(t)+W'_{l+1}(t)(x-t)+\dots\right)
\left(W_{l-1}(t)+W'_{l-1}(t)(x-t)+\dots\right)}
{\left(W'_{l}(t)(x-t)+\frac 12 W''_{l}(t)(x-t)^2+\dots\right)^2}=\\
& = &\frac{\left(W_{l+1}(t)+W'_{l+1}(t)(x-t)+\dots\right)
\left(W_{l-1}(t)+W'_{l-1}(t)(x-t)+\dots\right)}{W'^2_l(t)(x-t)^2}\times\\
&\times &\left(1-\frac{W''_l(t)}{W'_l(t)}(x-t)+\dots\right),
\end{eqnarray*}
where dots stand for the terms containing higher degrees of $(x-t)$.
Therefore the residue at $t$ is
$$
-\frac{W_{l-1}(t)W_{l+1}(t)W''_l(t)}{W'^3_l(t)}+
\frac{W'_{l+1}(t)W_{l-1}(t)+W'_{l-1}(t)W_{l+1}(t)}{W'^2_l(t)}\,,
$$
and the statement follows.  \hfill   $\triangleleft$

\medskip \noindent
{\bf Remark}\ \ Originally,  the equality (\ref{Ga}) was proved
 by A.~Gabrielov in a different, purely  algebraic, way (\cite{Ga}).

\section{The generating function of  a Schubert intersection}\label{s3}
\subsection{The Schubert intersection}\label{s31}
Consider the  intersection of Schubert varieties  $\I_{\{\w\}}(z)$ given by  (\ref{I}).
For any  $V\in \I_{\{\w\}}(z)$,  all finite singular points are
$\{ z_1,\dots, z_n\}$ and the Wronskian  is as in (\ref{W}). 
According to the definition of Schubert indices (\ref{E3}), 
the exponents $\rho_l(z_j)$ at $z_j$ and the exponents  $-d_l$ at $\infty$ 
are as follows,
\begin{equation}\label{E7}
\rho_l(z_j)=w_l(j)+p-l\,,\quad  d_l=d-w_l(n+1) +l-p\,,\quad 1\leq l\leq p,
\  1\leq j\leq n.
\end{equation}
We have
\begin{equation}\label{E5}
0\leq d_1<d_2<\dots <d_p=d,\quad
0=\rho_p(z_j)<\rho_{p-1}(z_j)<\dots <\rho_1(z_j)\leq d.
\end{equation}
The conditions $d_p=0$ and $\rho_p(z_j)=0$ for $1\leq j\leq n$ mean
that $V$ has no base point.
Denote $V_\bullet$  the  flag obtained by the  intersection of $V$
and  $\F_\bullet(\infty)$ (recall that we define  $\F_l(\infty)=\Poly_l$),
\begin{equation}\label{V}
V_\bullet\,=\,\left\{\,  
V_1\subset V_2 \subset \dots \subset V_p=V\,\right\},\quad  \dim V_l=l.
\end{equation}
The degrees of polynomials in $V_l$ are $d_1,\dots, d_l$.
Denote $W_l(x)$ the Wronskian of $V_l$, $1\leq l\leq p$.   
In particular,  $W_p(x)$ coincides with $W_{m,z}(x)$ given by  (\ref{W}). 

Define  polynomials
\begin{equation}\label{Z} 
Z_i(x)=\prod_{j=1}^n(x-z_j)^{m_j(i)},\quad 1\leq i\leq p-1, 
\end{equation}
where 
\begin{equation}\label{mij}
m_j(i) = \sum_{l=p+1-i}^p\rho_l(z_j)-\frac{i(i-1)}2\,=\,
\sum_{l=p+1-i}^{p}w_l(j),\ \ 1\leq i\leq p-1, \ 1\leq j\leq n\,.
\end{equation}
In particular, the condition $w_p(j)=0$ for $1\leq j\leq n$ gives
$m_j(1)=0$ and $Z_1(x)=1$.

\begin{lemma} 
The ratio $W_i(x)/Z_i(x)$ is a polynomial of degree
\begin{eqnarray}\label{k}
k_i& = &\sum_{l=1}^{p-i}d_l-\sum_{j=1}^n\sum_{l=i+1}^{p}\rho_l(z_j)
+(n-1)\frac{(p-i)(p-i-1)}2= \\
\ & = &i(d+1-p)-\sum_{l=1}^{i}w_l(n+1)-
\sum_{j=1}^{n}\sum_{l=p+1-i}^{p}w_l(j)\,, \  1\leq i\leq p-1.\nonumber
\end{eqnarray}
\end{lemma}

\noindent{\bf Proof:}\ \
For any $1\leq j\leq n$,  there exists a basis in $V_i$ consisting
of polynomials which have a root at $z_j$ of multiplicities at least
$\rho_p(z_j), \dots,  \rho_{p+1-i}(z_j)$. Therefore $W_i(x)$ has a root
at $z_j$ of order at least  $m_j(i)$.
Furthermore, $V_i$ is spanned by polynomials of degrees $d_1,\dots, d_i$,
therefore the degree of $W_i(x)$ is  $\sum_{l=1}^id_l-i(i-1)/2.$
Clearly $k_{p-i}=\deg W_i-m_1(i)-\dots -m_n(i)$.
\hfill   $\triangleleft$

\begin{cor}\label{cor} If the Schubert intersection $\I_{\{\w\}}(z)$
is non-empty, then the numbers $k_i$ given by {\rm (\ref{k})} are
nonnegative.
\end{cor}

\medskip
The definition of a nondegenerate $p$-plane, see Sec.~\ref{s13},
implies that if $V\in \I_{\{\w\}}(z)$ is nondegenerate, and if
$V_\bullet$ is the corresponding flag (\ref{V}),
then $V_l \in G_l(\Poly_{d_l})$,
the exponents of $V_l$ at $\infty$ of  are $-d_1,\dots, -d_l$, 
and the exponents at $z_j$ are
$\rho_p(z_j)=0, \rho_{p-1}(z_j),  \dots,  \rho_{p+1-l}(z_j)$,
for any $1\leq l\leq p$ and $1\leq j\leq n$. We arrive at
the following conclusion.

\begin{lemma}\label{Lemma}
Let $V$ be a nondegenerate $p$-plane in $\I_{\{\w\}}(z)$ and
$V_\bullet$ the flag {\rm (\ref{V})}.
Then  polynomials $T_{p-i}(x)= W_i(x)/Z_i(x)$, where
$W_i(x)$ is the Wronskian of $V_i$ and $Z_i(x)$ is defined by
{\rm (\ref{Z}), (\ref{mij})},  do not have multiple 
roots and satisfy
$$
T_{p-i}(z_j)\neq 0, \ \  1\leq j\leq n,\ \ 1\leq i\leq p-1; \ \ 
T_{p-1}(x)=W_1(x).
$$
Moreover, $T_{p-i}(x)$ and $T_{p-i-1}(x)$ do not have common roots
for $1\leq i\leq p-1$. \hfill   $\triangleleft$
\end{lemma}

\subsection{The generating function of  
$\I_{ \{\w\}}(z)$}\label{s32}
If $V$ is a nondegenerate $p$-plane in $ \I_{\{\w\}}(z)$, then 
Lemma~\ref{Lemma} implies that the relative discriminant of $W_i(x)$ 
with respect to $z$ is
$$
\Delta_{z}(W_i)=\Delta(T_{p-i}){\rm Res}^2(Z_i,T_{p-i}),\ \ 1\leq i\leq p-1,  
$$
and  the relative resultant   of  $W_i(x)$ and
 $W_{i+1}(x)$  with respect to $z$ is
$$
{\rm Res}_{z}(W_i, W_{i+1})={\rm Res}(T_{p-i}, T_{p-i-1})
{\rm Res}(T_{p-i}, Z_{i+1}){\rm Res}(T_{p-i-1}, Z_i),\ \  1\leq i\leq p-1,
$$
see Sec.~\ref{s13}.

\medskip
In order to re-write down the generating  function (\ref{Phi-w}) in terms of
unknown roots of polynomials
$T_1,\dots,T_{p-1}$, recall that if  $P(x)=(x-a_1)\dots (x-a_A)$ and
$Q(x)=(x-b_1)\dots(x-b_B)$,  then
$$
\Delta(P)=\prod_{1\leq i<j\leq A}(a_i-a_j)^2\,,\quad
{\rm Res}(P,Q)=\prod_{i=1}^A\prod_{j=1}^B(a_i-b_j).
$$

Denote
$$
t^{(i)}= \left(\, t_1^{(i)}\,,\  \dots\,,\  t_{k_i}^{(i)}\right)\,
\quad i=1,\dots, p-1,
$$
the roots of $T_i(x)$ and  write $\T=\left(\, t^{(1)}\,,\ \dots\,,
t^{(p-1)}\,\right)$, $\kk=(k_1,\dots, k_{p-1})$.
The numbers $m_j(i)$, $1\leq j\leq n$, $1\leq i\leq p-1$, are defined 
in (\ref{mij}).
Set  $m_j(0)=0$ and $m_j(p)=m_j$,  $1\leq j\leq n$, and  denote
 $\m=\{m_j(i),\ \ 1\leq j\leq n,\ \ 0\leq i\leq p\}$. We have

\begin{eqnarray}\label{Phi-t}
\Phi_{ \{\w\},z} (W_1,\dots, W_{p-1})\,&=&\, \Phi _{ \m,\kk, z}(\T)=\\
&=& \prod_{i=1}^{p-1}\ \
\prod_{1\leq l<s\leq k_i}\ \ (t^{(i)}_l-t^{(i)}_s)^2 \nonumber \\
&\times &\prod_{i=1}^{p-2}\ \prod_{l=1}^{k_{i}}\ \prod_{s=1}^{k_{i+1}}
(t^{(i)}_l-t_s^{(i+1)})^{-1} \nonumber\\
&\times &\prod_{i=1}^{p-1}\prod_{j=1}^n\ \prod_{l=1}^{k_i}\,
(t^{(i)}_l-z_j)^{2m_j(p-i)-m_j(p-i-1)-m_j(p-i+1)}\,.\nonumber
\end{eqnarray}

This is a rational function in $k_1+\dots+k_{p-1}$ variables.
If $k_i=0$ for  some $i$, then corresponding terms in (\ref{Phi-t}) are missing.

\medskip
According to (\ref{mij}), we get
\begin{equation}\label{exp}
2m_j(i)-m_j(i-1)-m_j(i+1)=\rho_{p+1-i}(z_j)-\rho_{p-i}(z_j)+1\leq 0,
\end{equation}
 for any $1\leq j\leq n$ and $1\leq i\leq p-1$.

The critical points with non-zero critical values are solutions to the equations
$$
\Phi^{-1}_{ \m,\kk, z}(\T)\ \frac{\partial\Phi_{\m,\kk,z}}{ \partial t_l^{(i)}}(\T)=0\,,
\quad 1\leq i\leq p-1, \quad 1\leq l\leq k_i\,.
$$
These equations have the following form,
\begin{eqnarray}\label{S}
\ &\ &\sum_{s\neq l}\frac 2{t_l^{(i)}-t_s^{(i)}}-
\sum_{s=1}^{k_{i-1}}\frac 1{t_l^{(i)}-t_s^{(i-1)}}-
\sum_{s=1}^{k_{i+1}}\frac 1{t_l^{(i)}-t_s^{(i+1)}}+\\
\ & + & \sum_{j=1}^{n}\frac {2m_j(p-i)-m_j(p-i-1)-m_j(p-i+1)}
{t_l^{(i)}-z_j}\,=\,0.\nonumber
\end{eqnarray}

The set of critical points is clearly invariant with respect to the group
$S^\kk=S^{k_1}\times\dots\times S^{k_{p-1}}\,$,
where  $S^{k_i}$ is the group of permutations of
$t_1^{(i)},\dots,t_{k_i}^{(i)}\,$. For any orbit,
the critical value is the same.

\subsection{Proof of Theorem 1}\label{s33}
Now we are in position to establish 
a one-to-one correspondence between the orbits
of critical points with non-zero critical value
of the function  $\Phi_{ \m,\kk, z}(\T)$ given by  {\rm (\ref{Phi-t})}
and the  nondegenerate
$p$-planes in  the intersection of Schubert varieties
$\I_{\{\w\}}(z)$ given by {\rm (\ref{I})}. The indices 
$\{\w\}$ and $\m,\kk$ are related as in (\ref{mij}) and (\ref{k}).

\medskip
 Let $V$ be a nondegenerate $p$-plane in
 $\I_{\{\w\}}(z)$.   Consider the
Wronskians $W_1,\dots, W_{p-1}$ corresponding to the
flag $V_\bullet$ defined by (\ref{V}), and write the equation $E_V$
in the form   (\ref{E4}).
Proposition~\ref{P3}  says that the roots
of polynomials $W_1,\dots, W_{p-1}$ which differ from the singular
points of the equation give a solution to the  system (\ref{S}).
Indeed, a direct calculation gives
\begin{eqnarray*}
\frac{W''_{p-i}(t_l^{(i)}) }{W'_{p-i}(t_l^{(i)}) } &=&
\sum_{s\neq l}\frac 2{t_l^{(i)}-t_s^{(i)}}
+ \sum_{j=1}^{n}\frac {2m_j(p-i)}{t_l^{(i)}-z_j}\,,\\
\frac{W'_{p-i\pm 1}(t_l^{(i)}) }{W_{p-i\pm 1}(t_l^{(i)}) }&=&
\sum_{s=1}^{k_{i\pm 1}}\frac 1{t_l^{(i)}-t_s^{(i\pm 1)}}
+\sum_{j=1}^{n}\frac {m_j(p-i\pm 1)}{t_l^{(i)}-z_j}\,.
\end{eqnarray*}

\medskip
Conversely,  any  critical point of $\Phi_{\m,\kk, z}(\T)$ with non-zero
critical value defines polynomials $T_{p-i}(x)$ and   
$W_i(x)=T_{p-i}(x)Z_i(x)$ with
$Z_i(x)$  given by (\ref{Z}), $1\leq i\leq p$.
The polynomials $W_i(x)$ determine an equation of the form (\ref{E4}).
Corollary~\ref{C3} gives $p$ linearly independent solutions 
$u_1,\dots, u_p$ to this equation,  and $W_i(x)$ is the Wronskian of
$u_1,\dots, u_p$. We have to show that  $\Span\{u_1,\dots, u_p\}$ is a  nondegenerate
$p$-plane in  the Schubert intersection $\I_{\{\w\}}(z)$ given by {\rm (\ref{I})}.

First, we check that $u_1,\dots, u_p$ are regular functions. 
For $u_1=W_1$ it is evident.
Fix  $1\leq i\leq p-1$ and consider the most interior integrand in the formula
for  $u_{i+1}$ of Corollary~\ref{C3}. It is $W_{i-1}W_{i+1}/W_i^2$ (recall that we set
$W_0=1$). This  is a rational function, hence it is regular at any point
which differs from the roots of $W_i$.
The roots of $W_i$ are  $z_j$ and $t_l^{(p-i)}$, $1\leq j\leq n$, 
$1\leq l\leq k_{p-i}$
According to (\ref{exp}),  $z_j$ is a root of the integrand
of  multiplicity $\rho_{p-i}(z_j)-\rho_{p-i+1}(z_j)-1\geq 0$,
and hence the integrand  is regular at $z_j$ for any $1\leq j\leq n$.
The residue at $t_l^{(p-i)}$  vanishes according to Proposition~\ref{P3}.
Therefore the integral has the form
$$
\int^{\xi}\frac{W_{i-1}(\tau)W_{i+1}(\tau)}{W^2_i(\tau)}d\tau
=\frac{c}{\left(\xi-t_l^{(p-i)}\right)}+ F(\xi),
$$
where $c$ is a constant number and   $F(\xi)$ is a function  regular at $t_l^{(p-i)}$.
Proposition~\ref{P3} says that $t_l^{(p-i)}$ is a simple root of  $W_i$,
hence the next integrand which is
$$
\frac{W_{i-2}(\xi)W_{i}(\xi)}{W_{i-1}^2(\xi)}
\left(\frac{c}{(\xi-t_l^{(p-i)})}+ F(\xi)\right)
$$
is regular at $t_l^{(p-i)}$.  We conclude that this integrand is regular at
all points which differ from the roots of  $W_{i-1}$, hence  we can repeat
the same arguments for this integrand as well.
In the last analysis,  we get that the function $u_{i+1}$ is regular, for
any $1\leq i\leq p-1$.

Next, the exponents of  the equation  with the solution space 
$V_i=\Span\{u_1,\dots, u_i\}$ 
at $z_j$ are exactly $\rho_l(z_j)$, $p-i+1\leq l\leq p$, $1\leq j\leq n$, and the
exponents at $\infty$ are $-d_1,\dots, -d_i$, where 
\begin{eqnarray}\label{dkm}
\rho_l(z_j)&= &m_j(p+1-l)-m_j(p-l) +p-l\,, \ \   j=1,\dots, n\,,\nonumber\\
d_l & = & k_{p-l}-k_{p+1-l}+l-1+\sum_{j=1}^n\left(m_j(l)-m_j(l-l)\right)\,,
 \ \ l=1,\dots p\,.
\end{eqnarray}
The relations (\ref{mij}) and (\ref{k}) assert that the exponents at the
points $z_1,\dots, z_n$ are non-negative integers, and the exponents at
$\infty$ are non-positive integers. Hence $u_1,\dots, u_p$
are polynomials and the solution space is a nondegenerate 
$p$-plane in $\I_{\{\w\}}(z)$.  \hfill   $\triangleleft$

\begin{cor}\label{C4}
The number of critical orbits with non-zero critical values
of the function $\Phi_{ \m,\kk, z}(\T)$ 
given by{\rm (\ref{Phi-t})} is at most the intersection number
$\sigma_{\w(1)}\cdot {\rm \ ... \ }\sigma_{\w(n)}
\cdot \sigma_{\w(n+1)}$, where $\{\w\}$ and $\m,\kk$ are related as in
{\rm (\ref{mij}),  (\ref{k})}.
\hfill   $\triangleleft$
\end{cor}

\subsection{Nondegenerate rational curves}\label{s34}
The  elements of   $G_p({\Poly}_d)$ can be 
interpreted as equivalence classes of rational curves. More precisely,
any basis $P_1,\dots, P_p$ of  $V\in G_p({\Poly}_d)$ defines
a {\it rational curve of degree} $d$ in $C\p^{p-1}$as follows,
$$
\gamma_V:\C\p^1\to\C\p^{p-1},\quad [t:s]\mapsto
\left[s^dP_1(t/s):\dots : s^dP_p(t/s)\right],
$$
and any other basis  defines {\it an equivalent curve}.
The singular points and the Schubert indices
(or, equivalently, the exponents) at these points
provide the {\it  ramification data} of the curve, see \cite{KhSo}.
 
The flag (\ref{V}) determine a set of $p-1$
{\it intermediate} curves $\gamma_{V(1)},\dots, \gamma_{V(p-1)}$;
the rational curve $\gamma_{V(l)}$ is defined by $V_l\in G_l({\Poly}_{d_l})$,
$1\leq l\leq p-1$.

If $\gamma_V$ is defined by a  nondegenerate $p$-plane $V$, then we say 
that it is a {\it nondegenerate} rational curve. A nondegenerate curve
with ramification data $z$ and $\{\w\}$ can be characterized as follows, 

\medskip
{\it for all intermediate curves, all ramification points
distinct from $z_1,\dots, z_n$ are flexes.}

\medskip
Theorem~1 can be reformulated in terms of rational curves.

\begin{cor}\label{C11}
There is a one-to-one correspondence between the critical points 
with non-zero critical values
of the function $\Phi_{\{\w\}, z}$ given by  {\rm (\ref{Phi-w})}
and the nondegenerate rational curves in $\C\p^{p-1}$ with the ramification data 
$z$ and $\{\w\}$.  \hfill   $\triangleleft$
\end{cor}

Our Conjecture says that for generic $z$, {\it all} rational functions
with prescribed ramification data  $z$ and $\{\w\}$ are nondegenerate.
The case $p=2$ has been done in \cite{S1}.

\section{Bethe vectors in the $\ls_p$ Gaudin model}\label{s4}

\subsection{Subspaces of singular vectors}(See \cite{FH}.)\label{s41}\ \
Consider the Lie algebra $\ls_p=\ls_p(\C)$. Denote  $(\cdot,\cdot)$  
the  Killing form on the dual to the Cartan subalgebra,
and  $\alpha_1,\dots, \alpha_{p-1}$ the simple roots. We have  
$$
(\alpha_i,\alpha_i)=2,\ \ (\alpha_i,\alpha_{i\pm 1})=-1,\ \
(\alpha_i,\alpha_j)=0 \ \ {\rm for\ \ } |i-j|\geq 2\,.
$$ 
Take a vector $\A=(a_1,\dots, a_{p-1})$ with nonnegative integer
coordinates, and denote $\Gamma_\A$
the irreducible  representation with integral dominant highest weight
$\Lambda_\A$,
$$(\Lambda_\A, \alpha_i)= a_i,\,\quad  1\leq i\leq p-1.$$
Consider the tensor product of $n$ irreducible representations
with integral dominant highest weights  $\Lambda_{\A(j)}$, 
$1\leq j\leq n$,
$$
 \Gamma_{\{\A\}}=\Gamma_{\A(1)}\otimes\dots\otimes \Gamma_{\A(n)}\,,
 \  \  \{\A\}=\{\A(1),\dots,\A(n)\}\,.
$$ 
Define the weight
$$
\Lambda(\kk)=\Lambda_{\A(1)}+\dots + \Lambda_{\A(n)}-
k_1\alpha_1-\dots-k_{p-1}\alpha_{p-1},
$$
where $\kk=(k_1,\dots, k_{p-1})$ is a vector with nonnegative
integer coordinates.
The subspace $\Sing_{\kk}\Gamma_{\{\A\}}$ of singular vectors
 of  the  weight  $\Lambda(\kk)$ in  $\Gamma_{\{\A\}}$ is as follows,
\begin{equation}\label{Sing}
\Sing_{\kk}\Gamma_{\{\A\}}=\{v\in \Gamma_{\{\A\}}\,\vert\,
h_iv=(\Lambda_\kk,\alpha_i)v,\ e_iv=0,
\ i=1,\dots,p-1\}\,,
\end{equation}
where  $\{e_i,f_i,h_i\}_{i=1}^{p-1}$ are the standard Chevalley  generators of $\ls_p$\,,
$$
[h_i,e_i]=2e_i\,,\ [h_i,f_i]=-2f_i\,,\ [e_i,f_i]=h_i\,;\
[h_i,h_j]=0\,,\ [e_i,f_j]=0\,  {\rm \ if \ } \ i\neq j.
$$
We will assume  $\Lambda(\kk)$ to be an integral dominant weight.
This means that the numbers
\begin{eqnarray}\label{Idw}
\left(\Lambda(\kk),\alpha_1\right)&=&\sum_{j=1}^n a_1(j)-2k_1+k_2,\\
\left(\Lambda(\kk),\alpha_l\right)&=&
\sum_{j=1}^n a_l(j)+k_{l-1}-2k_{l}+k_{l+1},
\quad  2\leq l \leq p-2, \nonumber\\
\left(\Lambda(\kk),\alpha_{p-1}\right)&=&
\sum_{j=1}^n a_{p-1}(j)+k_{p-2}-2k_{p-1}\nonumber
\end{eqnarray}
are nonnegative integers.
Denote $\A(n+1)$ the vector with coordinates
$$
a_l(n+1)=\left(\Lambda(\kk),\alpha_l\right),\quad 1\leq l\leq p-1,
$$
and $\Gamma_{\A(n+1)}$ the $\ls_p$-representation with highest
weight $\Lambda(\kk)$.
The dimension of $\Sing_{\kk}\Gamma_{\{\A\}}$ is
the multiplicity of $\Gamma_{\A(n+1)}$
in the decomposition of $\Gamma_{\{\A\}}$
into the direct sum of irreducible representations.

\subsection{The master function of the  Gaudin model associated
with $z$ and $\Sing_{\kk}\Gamma_{\{\A\}}$}\label{s42}
The master function determines those of common eigenvectors of 
the Gaudin hamiltonians  (\ref{H}) which are singular vectors of  
the weight $\Lambda(\kk)$ in $\Gamma_{\{\A\}}$.
This is a function in $k_1+\dots+k_{p-1}$ complex variables
$$
\{t_l^{(i)},\quad 1\leq i\leq p-1,\quad 1\leq l\leq k_i\}
$$
which has the following form (see \cite{ReV}),
$$
\prod_{i=1}^{p-1}\prod_{l=1}^{k_i}\prod_{j=1}^n
(t_l^{(i)}-z_j)^{-a_i(j)}
\prod_{i=1}^{p-1}\prod_{1\leq l<s\leq k_i}
(t_l^{(i)}-t_s^{(i)})^2
\prod_{i=1}^{p-2}\prod_{l=1}^{k_{i}}
\prod_{s=1}^{k_{i+1}}(t_l^{(i)}-t_s^{(i+1)})^{-1}.
$$
We see that this is the same function as in  (\ref{Phi-t})  with
$$
a_i(j)=m_j(p-i-1)-2m_j(p-i)+m_j(p-i+1),\quad 1\leq i\leq p-1,\ \ 1\leq j\leq n.
$$
Consider the Schubert intersection $\I_{\{\w\}}(z)$ given by (\ref{I}).
The indices $\w(j)$, $1\leq j\leq n+1,$ determine
numbers $k_i$ and $m_j(i)$, $1\leq i\leq p-1$, $1\leq j\leq n$, in accordance
with (\ref{mij}) and (\ref{k}). We arrive at the following result.

\begin{cor}\label{C31} There is a one-to-one correspondence between
the nondegenerate $p$-planes in $\I_{\{\w\}}(z)\subset G_p({\Poly}_d)$ and
the Bethe vectors of the Gaudin model associated with $z$ and 
$\Sing_{\kk}\Gamma_{\{\A\}}$,
where
\begin{eqnarray*}
a_i(j)&=&w_i(j)-w_{i+1}(j),\quad
1\leq i\leq p-1,\quad 1\leq j\leq n;\\
k_i&=&i(d+1-p)-\sum_{l=1}^{i}w_l(n+1)-
\sum_{j=1}^{n}\sum_{l=p+1-i}^{p}w_l(j)\,,\ \   1\leq i\leq p-1.
\quad \hfill \triangleleft
\end{eqnarray*}
\end{cor}

\medskip
For all $1\leq j\leq n$,
the Schubert class $\sigma_{\w(j)}$ and the $\ls_p$-module
$\Gamma_{\A(j)}$ correspond to the same  Young diagram, see \cite{F}.
This diagram has $p-1$ rows with $w_i(j)$ boxes in  the $i$-th row,
or equivalently, it has $a_i(j)$ columns of the length $i$.

\medskip
As one can easily check, the condition that $\Lambda(\kk)$ given by
(\ref{Idw}) is an integral dominant weight is necessary for 
$\I_{\{\w\}}(z)$ is non-empty, see Corollary~\ref{cor}.

\medskip
Recall that the Schubert classes $\sigma_{\tilde\w}$
and  $\sigma_{\w}$ in $G_p({\Poly}_d)$ are
{\it dual} if  $|\tilde\w|+|\w|=p(d+1-p)$ and
the intersection number $\sigma_{\w}\cdot\sigma_{\tilde\w}$
is $1$ in  $H^{2p(d+1-p)}(G_p({\Poly}_d))$.
The Schubert class dual to $\sigma_{\w(n+1)}$ and  $\Gamma_{\A(n+1)}$
correspond to the same  Young diagram, \cite{F}.
The relation of the Schubert calculus in  to
the representation theory of Lie algebra $\ls_p(\C)$
implies the following  fact.

\begin{prop}\label{P5}  {\rm \cite{F}}\ \ 
The multiplicity of  $\Gamma_{\A(n+1)}$
in $\Gamma_{\{\A\}}$ is the intersection number of the Schubert classes
$\sigma_{\w(1)}\cdot {\rm \ ... \ }\sigma_{\w(n)}
\cdot \sigma_{\w(n+1)}$.  \hfill  $\triangleleft$
\end{prop}

This Proposition and Corollary~\ref{C4} give the Corollary 
from Theorem~1 of Sec.~\ref{s14}.

\subsection{Eigenvalues of Bethe vectors}\label{s43}
Our construction provides a
correspondence between Bethe vectors and Fuchsian differential 
equations with only polynomial solutions. 
On the other hand, as it is explained in Sec.~5.4--5.8 of~\cite{Fr}, the eigenvalues
$\mu=(\mu_1,\dots, \mu_n)$ of any common eigenvector of the Gaudin
hamiltonians  (\ref{H}) determine an $\ls_p$-oper which
generates trivial monodromy representation of the fundamental group
$\pi_0\left(\C\p^1\setminus\{z_1,\dots, z_n,\infty\}\right)$.
This oper is in fact a differential operator of the form
$$
\D_\mu = \left(\frac d{dx}\right)^p+f_{p-2}(x)\left(\frac d{dx}\right)^{p-2}
+\dots +f_0(x)\,
$$
with meromorphic coefficients. 

The  operator $\D_\mu$ corresponds to  a certain Fuchsian
differential equation with only polynomial solutions and  the 
singular points belonging to the set $\{z_1,\dots, z_n,\infty\}$.
More precisely, the link is as follows,
$$
{\rm Ker\,}\D_\mu=\left\{\, \frac{P(x)}{\left({W_V(x)}\right)^{1/p}}\,, \quad
P(x)\in V\, \right\}\,,
$$
where  $V$ is a certain $p$-dimensional subspace in the vector space of complex
polynomials such that  all finite singular points of $V$ belong to the set
$\{z_1,\dots, z_n\}$.

In  Sec.~5 of~\cite{Fr}, it is proved 
that the existence of a common eigenvector of the Gaudin
hamiltonians with eigenvalues $\mu=(\mu_1,\dots, \mu_n)$
implies the existence of a Bethe vector with these eigenvalues.
Presumably, the Bethe vectors of the Gaudin model span the relevant
subspace of the singular vectors. For the case of $\ls_2$, this
has been  proved in \cite{ReV}. The simplicity of the spectrum
of the Gaudin hamiltonians then can be easily deduced from the relation
to differential equations, see \cite{S2}. Likely, the connection between 
the Frenkel construction and the presented one implies that the Bethe vectors 
are separated by their eigenvalues. Then the simplicity of the spectrum is equivalent
to the statement that the Bethe vectors span the subspace of singular vectors.

\section{Bethe vectors in the tensor product of symmetric powers of
the standard $\ls_p$-representation}\label{s5}

\subsection{Special Schubert varieties and Fuchsian equations}\label{s51}
Consider  the Schubert intersections such that the Schubert indices corresponding 
to all finite singular points are special, that is $\w(j)=(m_j,0,\dots,0)$ for certain
positive integers $m_j$, $1\leq j\leq n$. The Young diagram of the special
Schubert variety  $\Omega_{m_j}$ consists of one row with $m_j$ boxes.

\medskip
Denote  $\dd =(d_1,\dots, d_p)$, where integers $d_1,\dots, d_p$ satisfy
\begin{equation}\label{dd}
0\leq d_1<\dots<d_{p-1}<d_p=d,\quad
d_1+\dots+ d_p=m_1+\dots +m_n + \frac {p(p-1)}2.
\end{equation}
The Schubert index at infinity $\w(n+1)$ is defined by $\dd$. 
We will denote $\I_{\{m,\dd\}}(z)$ the intersection (\ref{Im}) with
$w_l(n+1)=d-d_l+l-p$,  $1\leq l\leq p$,  see  (\ref{E3}).

\begin{prop}\label{P7}
If $V\in \I_{\{m,\dd\}}(z)$, then the equation $E_V$ has the form
\begin{equation}\label{E8}
\prod_{j=1}^n(x-z_j)u^{(p)}(x)+F_{1}(x)u^{(p-1)}(x)+\dots +F_p(x)u(x)=0,
\end{equation}
where $F_i(x)$ is a polynomial of degree at most $n-i$ for any
$1\leq i\leq p$, and
\begin{equation}\label{E9}
F_1(x)=\left(-\frac{m_1}{x-z_1}-\dots-
\frac{m_n}{x-z_n}\,\right)\cdot \prod_{j=1}^n(x-z_j)\,.
\end{equation}
Conversely, if $m_1,\dots, m_n$ are positive
integers,  if all solutions to the equation
{\rm (\ref{E8}),  (\ref{E9})} are polynomials,
and if the degrees $d_1,\dots, d_p$ of solutions satisfy {\rm (\ref{dd})},
then the solution space belongs to  $ \I_{\{m,\dd\}}(z)$.
\end{prop}

\medskip\noindent{\bf Proof :}\ \
The decomposition of the determinant in the left hand side of
(\ref{E1}) with respect to the first column gives the following form of the
equation $E_V$,
\begin{eqnarray*}
&&W_V(x)u^{(p)}(x)-W'_V(x)u^{(p-1)}(x)+F_{p-2}(x)u^{(p-2)}(x)+\\
&&+\dots +F_{2}(x)u''(x)+F_{1}(x)u'(x)+F_0(x)u(x)=0\,,\nonumber
\end{eqnarray*}
where $F_0,\dots, F_{p-2}$ are suitable polynomials.
We see  that the multiplicities
of $z_j$ as a root of  the coefficient at $u^{(p)}$ and of
the coefficient at $u^{(p-1)}$ differ by 1. Furthermore,
for any $1\leq j\leq n$ we have $\rho_l(z_j)=p-l,\ 2\leq l\leq p$,
in accordance with the definition of the Schubert index (\ref{E3}).
Therefore it is enough to proof that if the equation
$$
x^kF_0(x)u^{(p)}(x)+x^{k-1}F_1(x)u^{(p-1)}(x)+
F_2(x)u^{(p-2)}(x)+\dots+F_p(x)u(x)=0,
$$             
where $F_0(0)\neq 0$,  $F_1(0)\neq 0$, and $k\geq 2$,
has solutions of the form
$$
x^if_i(x),\quad f_i(0)\neq 0,\quad i=0,1,\dots, p-2,
$$
then $F_2(0)=\dots=F_p(0)=0$. Indeed, the
substitution of $u(x)=x^{p-2}f_{p-2}(x)$ into the equation
gives $F_2(0)=0$, then the substitution of $u(x)=x^{p-3}f_{p-3}(x)$ into
the equation leads to  $F_3(0)=0$ and so on. Finally, the substitution
of $u(x)=f_0(x)$ gives $F_p(0)=0$ and finishes the proof of the first part
of the Proposition.

The second part follows from Proposition~\ref{P1}, since the exponents
of the equation  \rm (\ref{E8}),  (\ref{E9}) at all singular points are known.
\hfill $\triangleleft$

\medskip
As we explained in Sec.~\ref{s15}, the generating function of  the intersection
$\I_{\{m,\dd\}}(z)$ has the form (\ref{G}). 
According to Lemma~\ref{Lemma}, we have $W_i(x)=T_{p-i}(x)$ for
$1\leq i\leq p-1$. Thus in terms of  unknown  roots of  the Wronskians 
$W_1,\dots, W_{p-1}$, we get
\begin{eqnarray}\label{P}
\Phi (\T)=\Phi_{\kk,m,z}(\T)&= &\prod_{j=1}^n\ \prod_{l=1}^{k_{1}}\,
(t^{(1)}_l-z_j)^{-m_j}\\
&\times&\prod_{i=1}^{p-1}\ \prod_{1\leq l<s\leq k_l}\ \
(t^{(i)}_l-t^{(i)}_s)^2 \nonumber\\
&\times &\prod_{i=1}^{p-2}\ \prod_{l=1}^{k_{i}}\ \prod_{s=1}^{k_{i+1}}
(t^{(i)}_l-t_s^{(i+1)})^{-1} \nonumber \,,
\end{eqnarray}
here  $t^{(i)}=(t^{(i)}_1,\dots, t^{(i)}_{k_i})$ are the roots of $W_{p-i}(x)$, 
$1\leq i\leq p-1$, in accordance with  notation of Sec.~\ref{s32}.

\subsection{Subspaces of singular vectors in the tensor product
of symmetric powers of the standard $\ls_p$-representation}\label{s52}

Denote $\LL$ the standard $\ls_p$-representation.
Fix vector $m=(m_1,\dots,m_n)$  with positive integer coordinates
and set
$$
L_j=\Sym^{m_j}\LL,\quad  \LL^{\otimes m}=L_1\otimes\dots \otimes L_n\,.
$$
With notation of Sec.~\ref{s41}, we have $L_j=\Gamma_{(m_j,0,\dots,0)}$.
The corresponding Young diagram is the same as for the special
Schubert variety $\Omega_{m_j}$;  it is a row with  $m_j$ boxes.
Denote $E=(1,0,\dots,0)$. The  weight of $L_j$ is $m_jE $. 
For nonnegative integers $k_1,\dots,k_{p-1}$,
consider the weight
$$\Lambda=(m_1+\dots+m_n)E-k_1\alpha_1-\dots k_{p-1}\alpha_{p-1}\,,
\ \ \kk=(k_1,\dots,k_{p-1})\,.
$$
Set $k_0=m_1+\dots+m_n$. We have 
$$
(\Lambda,\alpha_i)=k_{i-1}-2k_{i}+k_{i+1}\,,  \quad 1\leq i\leq p-2,
\quad (\Lambda,\alpha_{p-1})=k_{p-2}-2k_{p-1}.
$$
Clearly $\Lambda$ is an integral dominant weight if and
only if
$$
k_{i-1}+k_{i+1}\geq 2k_{i}\,,\quad 1\leq i\leq p-2,
\quad k_{p-2} \geq 2k_{p-1}\,.
$$
On the other hand, the relation (\ref{dkm}) in our case takes the form
$$
d_1=k_{p-1}\,,\quad d_i=k_{p-i}-k_{p+1-i}+i-1\,,\quad 1\leq l\leq p,
$$
as for any $1\leq j\leq n$ we have $m_j(i)=0$, $1\leq i\leq p-1$.
Thus  $\Lambda$ is an integral dominant
weight if and only if  $\dd=(d_1,\dots,d_p)$ satisfies (\ref{dd}).

\medskip
Consider  $\Sing_{\kk}\LL^{\otimes m}$,
the subspace of singular  vectors  of the weight $\Lambda$ in  $\LL^{\otimes m}$.
The master function of  the Gaudin model associated with $z$ and  
$\Sing_{\kk}\LL^{\otimes m}$ is given by  (\ref{P}).
If $k_i$ vanishes for some $i$,
then the corresponding terms in $\Phi(\T)$ are missing.

\medskip
According to  Corollary~\ref{C4} and Proposition~\ref{P5}, 
the number of critical orbits of the function $\Phi_{\kk,m,z}(\T)$
is bounded from above by  the dimension of $\Sing_{\kk}\LL^{\otimes m}$.

\subsection{Case $n=2$}\label{n=2} In this case we can suppose
  $z=(0,1)$, and the special Schubert intersection takes the form
\begin{eqnarray*}
&&\I_{\{(m_1,m_2),\dd\}}(0,1)=\Omega_{m_1}(0)\cap\Omega_{m_2}(1)
\cap\Omega_{\w}(\infty)\subset G_p({\Poly}_d)\,,\\
&&\w=(w_1,\dots,w_p),\quad \dd=(d_1,\dots,d_p), \quad w_l=d-d_l+l-p.
\end{eqnarray*}
\begin{prop}\label{PC5}
 $\I_{\{(m_1,m_2),\dd\}}(0,1)$ consists of a single 
 element if $\dd=(0,1,\dots, p-3,  m_1+m_2-d+2p-3,d)$
and is empty otherwise. If  $V\in  \I_{\{(m_1,m_2),\dd\}}(0,1)$,
then the equation $E_V$  has the form
\begin{eqnarray*}
&\ & x(x-1)u^{(p)}(x)+\left((-m_1-m_2)x+m_1\right)u^{(p-1)}(x)+\\
&+& (d-p+2)(m_1+m_2+p-d-1)u^{(p-2)}(x)=0.
\end{eqnarray*}    
\end{prop}

\noindent{\bf Proof:}\ \ 
In our case the equation of Proposition~\ref{P7} 
takes the form
$$
x(x-1)u^{(p)}(x)+\left((-m_1-m_2)x+m_1\right)u^{(p-1)}(x)+
cu^{(p-2)}(x)=0\,,
$$
where $c$ is some constant. Clearly  the functions $1,\, x,\dots, x^{p-3}$ 
are solutions to this equation, i.e.  $d_l=l-1$, $1\leq l\leq p-2$.  The  value of  $d_{p-1}$
is defined  then by  (\ref{dd}),  $d_{p-1}= m_1+m_2-d+2p-3$. We get 
$$\w=(d+1-p, \dots, d+1-p, 2(d+1-p)-m_1-m_2,0).$$
It remains to notice that the obtained equation is the Gauss hypergeometric equation
with respect to $u^{(p-2)}(x)$,  \cite{R}. Therefore $c$ is the product of the 
corresponding exponents at infinity, which are $-d_{p-1}+(p-2)$ and $-d+(p-2)$.
\hfill $\triangleleft$

\medskip
Denote  $\Phi^0(\T)$ the function given by
(\ref{P}) with $n=2$, $m=(m_1,m_2)$,  $(z_1,z_2)=(0,1)$.
This is the master function of  the Gaudin model associated 
with $(0,1)$ and the subspace of singular vectors of the weight
$$
\Lambda_0=(m_1+m_2)E-k_{1}\alpha_1-\dots - k_{p-1}\alpha_{p-1}
$$
in  $\Sym^{m_1}\LL\otimes\Sym^{m_2}\LL$. 
The Pieri formula in this case can be formulated as follows,  \cite{F}.

\begin{lemma}\label{L7}
The dimension of the subspace of singular vectors of the weight $\Lambda_0$ 
in $\Sym^{m_1}\LL\otimes\Sym^{m_2}\LL$
is $1$ if $0\leq k_{1}\leq \min\{m_1,m_2\}$ and $k_2=\dots=k_{p-1}=0$,
and $0$ otherwise. \hfill   $\triangleleft$
\end{lemma}

\begin{prop}\label{L8}
 The number of  critical orbits with non-zero
critical values of the function $\Phi^0(\T)$ is $1$ if  the highest weight $\Lambda_0$
enters the tensor product $\Sym^{m_1}\LL\otimes\Sym^{m_2}\LL$, and $0$ otherwise.
\end{prop}

\noindent{\bf Proof:}\ \ 
If $V$ is an element of $G_p({\Poly}_d)$ corresponding to a critical point
of $\Phi^0(\T)$, then  Proposition \ref{PC5} asserts that one can take
$1,\, x,\dots, x^{p-3}$ as the first $p-2$  polynomials of the basis of $V$.
Thus  the Wronskians of the first $p-2$ subspaces of $V_\bullet$, see (\ref{V}), 
are  $W_1(x)=\dots =W_{p-2}(x)=1$, i.e. $k_2=\dots=k_{p-1}=0$.   

For  $k_2=\dots=k_{p-1}=0$, the function  $\Phi^0(\T)$  
coincides  with the master function of the $\ls_2$ Gaudin model associated with 
$z=(0,1)$ and the subspace of singular vectors of the weight $m_1+m_2-2k_{1}\geq 0$ 
in the tensor product of  two irreducible $\ls_2$-representations with highest weights 
$m_1$ and $m_2$.  Results of   Sec.~9  of~\cite{ReV} say that critical points with 
non-zero critical value exist only if $0\leq k_{1}\leq \min\{m_1,m_2\}$;  moreover
all of them are nondegenerate, lie in the same orbit and define a non-zero Bethe vector.
\hfill   $\triangleleft$

\begin{cor}
If  $\I_{\{(m_1,m_2),\dd\}}(0,1)$ is non-empty, then it is a nondegenerate $p$-plane.
 \hfill   $\triangleleft$
 \end{cor}

\subsection{Proof of  Theorem~2}\label{s53}
We showed that the  number of critical orbits with non-zero critical values
is $\dim \Sing_{\kk}\LL^{\otimes m}$ 
in the case  $m=(m_1,m_2)$ and $z=(0, 1)$.

Results of N.~Reshetikhin and A.~Varchenko
(see Theorem~9.16  and Theorem~10.4 of \cite{ReV})  imply that 
 $\dim \Sing_{\kk}\LL^{\otimes m}$  gives then a  bound from below for the
 number of critical orbits with non-zero critical values  of the function  
 $\Phi (\T)=\Phi_{\kk,m,z}(\T)$ given by
 (\ref{P}),  for generic $z=(z_1,\dots, z_n)$.
 Taking into account the Corollary from Theorem~1,  we get that the number of  
 these critical orbits is  exactly  $\dim \Sing_{\kk}\LL^{\otimes m}$.
Moreover,   arguments similar to those  in the proof of  Theorem~9.16 in~\cite{ReV}
show that distinct orbits define distinct non-zero Bethe vectors.
This finishes the proof of Theorem~2.  \hfill   $\triangleleft$

\begin{cor}\label{C9}
For generic $z=(z_1,\dots,z_n)$, the number of equations of the form
{\rm (\ref{E8}),  (\ref{E9})}
having polynomial  solutions of degrees $d_1,\dots,d_p$ with
$0\leq d_1 <\dots < d_p$ equals the dimension of
$\Sing_{\kk}\LL^{\otimes m}$, where $\kk=(k_1,\dots,k_{p-1})$, 
$$
k_i=d_1+\dots+d_{p-i}+(n-1)\frac{(p-i)(p-i-1)}2\,,
\quad i=1,\dots, p-1.\quad \hfill   \triangleleft
$$
\end{cor}

\newpage

\  \

\end{document}